\documentclass[11pt]{amsart}

\renewcommand{\bar}{\overline}

\newcommand{\pa}{\partial}
\renewcommand{\phi}{\varphi}

\newcommand{\pz}{\partial_z}
\newcommand{\pzb}{\partial_{\bar z}}

\usepackage{amsmath,amsthm,amscd}
\usepackage{calc}

\title
[]{Geometric Aspects of the Moduli Space of Riemann Surfaces}
\author{Kefeng Liu}
\address{Department of Mathematics\\
University of California at Los Angeles\\ Los Angeles, CA 90095-1555, USA\\
Center of Mathematical Sciences, Zhejiang University, Hangzhou,
China} \email{liu@math.ucla.edu, liu@cms.zju.edu.cn}

\author{Xiaofeng Sun}
\address{Department of Mathematics\\
Harvard University\\ Cambridge, MA 02138, USA}
\email{xsun@math.harvard.edu}
\author{Shing-Tung Yau}
\address{Department of Mathematics\\
Harvard University\\ Cambridge, MA 02138, USA}
\email{yau@math.harvard.edu}
\date{\today}

\thanks{The authors are supported by the NSF}

\newtheorem{theorem}{Theorem}[section]

\newtheorem{lemma}{Lemma}[section]
\newtheorem{cor}{Corollary}[section]
\newtheorem{prop}{Proposition}[section]

\newtheorem{definition}{Definition}[section]

\theoremstyle{remark}
\newtheorem{rem}{Remark}[section]

\begin{document}
\maketitle

\numberwithin{equation}{section}

\tableofcontents

\newcommand{\M}{{\mathcal M}}
\section{Introduction}\label{intro}

The study of moduli space and Teichm\"uller space has a long
history. These two spaces lie in the intersections of the researches
in many areas of mathematics and physics. Many deep results have
been obtained in history by many famous mathematicians. Here we will
only mention a few that are closely related to our discussions.

 Riemann was the first who considered the space $\M$ of all complex
structures on an orientable surface modulo the action of orientation
preserving diffeomorphisms. He derived the dimension of this space
\[
\dim_{\mathbb R}\M=6g-6
\]
where $g\geq 2$ is the genus of the topological surface.

In 1940's, Teichm\"uller considered a cover of $\M$ by taking the
quotient of all complex structures by those orientation preserving
diffeomorphims which are isotopic to the identity map. The
Teichm\"uller space $\mathcal T_g$ is a contractible set in $\mathbb
C^{3g-3}$. Furthermore, it is a pseudoconvex domain. Teichm\"uller
also introduced the Teichm\"uller metric by first taking the $L^1$
norm on the cotangent space of $\mathcal T_g$ and then taking the
dual norm on the tangent space. This is a Finsler metric. Two other
interesting Finsler metrics are the Carath\'eodory metric and the
Kobayashi metric. These Finsler metrics have been powerful tools to
study the hyperbolic property of the moduli and the Teichm\"uller
spaces and the mapping class groups. For example in 1970's Royden
proved that the Teichm\"uller metric and the Kobayashi metric are
the same, and as a corollary he proved the famous result that the
holomorphic automorphism group of the Teichm\"uller space is exactly
the mapping class group.

Based on the Petersson pairing on the spaces of automorphic forms,
Weil introduced the first Hermitian metric on the Teichm\"uller
space, the Weil-Petersson metric. It was shown by Ahlfors that the
Weil-Petersson metric is K\"ahler and its holomorphic sectional
curvature is negative. The works of Ahlfors and Bers on the
solutions of Beltrami equation put a solid fundation of the theory
of Teichm\"uller space and moduli space \cite{ab1}. Wolpert studied
in details the Weil-Petersson metric including the precise upper
bound of its Ricci and holomorphic sectional curvature. From these
one can derive interesting applications in algebraic geometry. For
example, see \cite{liu1} .

Moduli spaces of Riemann surfaces have also been studied in details
in algebraic geometry since 1960. The major tool is the geometric
invariant theory developed by Mumford. In 1970's, Deligne and
Mumford studied the projective property of the moduli space and they
showed that the moduli space is quasi-projective and can be
compactified naturally by adding in the stable nodal surfaces
\cite{dm1}. Fundamental works have been done by Gieseker, Harris and
many other algebraic geometers.

The work of Cheng-Yau \cite{cy2} in the early 80s showed that
there is a unique complete K\"ahler-Einstein metric on the
Teichm\"uller space and is invariant under the moduli group action.
Thus it descends to the moduli space. As it is well-known, the
existence of the K\"ahler-Einstein metric gives deep algebraic geometric
results, so it is natural to understand its properties like the
curvature and the behaviors near the compactification divisor. In
the early 80s, Yau conjectured that the K\"ahler-Einstein metric is
equivalent to the Teichm\"uller metric and the Bergman metric
\cite{cy2}, \cite{yau2}, \cite{yau4}.

In 2000, McMullen introduced a new metric, the McMullen metric by
perturbing the Weil-Petersson metric to get a complete K\"ahler
metric which is complete and K\"ahler hyperbolic. Thus the lowest
eigenvalue of the Laplace operator is positive and the
$L^2$-cohomology is trivial except for the middle dimension
\cite{mc}.

The moduli space appears in many subjects of mathematics, from
geometry, topology, algebraic geometry to number theory. For
example, Faltings' proof of the Mordell conjecture depends heavily
on the moduli space which can be defined over the integer ring.
Moduli space also appears in many areas of theoretical physics. In
string theory, many computations of path integrals are reduced to
integrals of Chern classes on the moduli space. Based on conjectural
physical theories, physicists have made several amazing conjectures
about generating series of Hodge integrals for all genera and all
marked points on the moduli spaces. The proofs of these conjectures
supply strong evidences to their theories.

 Our goal of this project is to understand the
geometry of the moduli spaces. More precisely, we want to
understand the relationships among all of the known canonical
complete metrics introduced in history on the moduli and the
Teichm\"uller spaces, and more importantly to introduce new
complete K\"ahler metrics with good curvature properties: the
Ricci metric and the perturbed Ricci metric. Through a detailed
study we proved that these new metrics have very good curvature
properties and very nice Poincar\'e-type asymptotic behaviors
\cite{lsy1}, \cite{lsy2}. In particular we proved that the
perturbed Ricci metric has bounded negative Ricci and holomorphic
sectional curvature and has bounded geometry. To the knowledge of
the authors this is the first known such metric on moduli space
and the Teichm\"uller spaces with such good properties. We know
that the Weil-Petersson metric has negative Ricci and holomorphic
sectional curvature, but it is incomplete and its curvatures are
not bounded from below. Also note that one has no control on the
signs of the curvatures of the other complete K\"ahler metrics
mentioned above.

 We have obtained a series of results. In \cite{lsy1} and \cite{lsy2} we have
proved that all of these known complete metrics are actually
equivalent, as consequences we proved two old conjectures of Yau
about the equivalence between the K\"ahler-Einstein metric and the
Teichm\"uller metric and also its equivalence with the Bergman
metric. In both \cite{yau2} and \cite{yau4} which were both
written in early 80s, Yau raised various questions about the
K\"ahler-Einstein metric on the Teichm\"uller space. By using the
curvature properties of these new metrics, we obtained good
understanding of the K\"ahler-Einstein metric such as its boundary
behavior and the strongly bounded geometry. As one consequence we
proved the stability of the logarithmic extension of the cotangent
bundle of the moduli space \cite{lsy2}. Note that the major parts
of our papers were to understand the K\"ahler-Einstein metrics and
the two new metrics. One of our goal is to find a good metric with
the best possible curvature property. The perturbed Ricci metric
is close to be such metric. We hope to understand its Riemannian
curvature in the future. The most difficult part of our results is
the study of the curvature properties and the asymptotic behaviors
of the new metrics near the boundary, only from which we can
derive geometric applications such as the stability of the
logarithmic cotangent bundle. The comparisons of those classical
metrics as well as the two new metrics are quite easy and actually
simple corollaries of the study and the basic definitions of those
metrics. In particular the argument we used to prove the
equivalences of the Bergman metric, the Kobayashi metric and the
Carath\'eodory metric is rather simple from basic definitions and
Yau's Schwarz lemma, and is independent of the other parts of our
works.

Our first paper was post in the webpage since February 2004 and
widely circulated. Since then the first and the second author have
given several lectures about the main results and key ideas of
both of our papers. In July we learned from the announcement of
S.-K. Yeung in Hong Kong University
\footnote{http://hkumath.hku.hk/$\sim$imr/records0304/GEO-YeungSK.pdf}
 where he announced he could prove a small and easy part of our results about
the equivalences of some of these metrics by using a bounded
pluri-subharmonic function. \footnote{ We received a hard copy of
Yeung's paper in November 2004 where he used a similar method to
ours in \cite{lsy2} to compare the Bergman, the Kobayashi and the
Carath\'eodory metric. It should be interesting to see how one can
use the bounded psh function to derive these equivalences.}

The purpose of this note is to give a brief overview of our
results and their background. It is based on the lecture delivered
by the first author in the First International Conference of
Several Complex Variables held in the Capital Normal University in
August 23-28, 2004. All of the main results mentioned here are
contained in \cite{lsy1} and \cite{lsy2} which the interested
reader may read for details. They have been circulated for a
while. The first author would like to thank the organizers for
their invitation and hospitality.

\section{The Topological Aspects of the Moduli Space}
The topology of the Teichm\"uller space is trivial, since it is
topologically a ball. But how to compactify it in a natural and
useful way is still an interesting problem. Penner has done
important works on this problem. The compactification of
Teichm\"uller space is useful in three dimensional topology. The
topology of the moduli space and its compactification is highly
nontrivial and have been well-studied for the past years from many
point of views. Here we only mention the recently proved Mumford
conjecture about the stable cohomology of the moduli spaces; the
Witten conjecture about the KdV equations for
 the generating series of the integrals of the $\psi$ classes; the
Mari\~no-Vafa conjecture about the closed expressions for the
generating series of triple Hodge integrals.

The first two results mentioned above are already well-known. Here
we would like to explain a little more details about the
Mari\~no-Vafa conjecture
 proved in \cite{llz} which gives a closed formula
 for the generating series of triple Hodge
integrals of all genera and all possible marked points, in terms of
Chern-Simons knot invariants.

Hodge integrals are defined as the intersection numbers of $\lambda$
classes and $\psi$ classes on the Deligne-Mumford moduli spaces of
stable Riemann surfaces $\bar\M_{g,h}$, the moduli with $h$ marked
points. Recall that a point in $\bar\M_{g,h}$ consists of
$(C,x_1,\ldots,x_h)$, a (nodal) Riemann surface $C$ and $h$ smooth
points on $C$.

The Hodge bundle $\mathbb{E}$ is a rank $g$ vector bundle over
$\bar\M_{g,h}$ whose fiber over $[(C,x_1,\ldots,x_h)]$ is
$H^0(C,\omega_C)$. The $\lambda$ classes are the Chern Classes:
$$
\lambda_i=c_i(\mathbb{E})\in H^{2i}(\bar\M_{g,h};\mathbb Q).
$$

On the other hand the cotangent line $T_{x_i}^* C$ of $C$ at the
$i$-th marked point $x_i$ gives a line bundle $\mathbb{L}_i$ over
$\bar\M_{g,h}$. The $\psi$ classes are also Chern classes:
$$
\psi_i=c_1(\mathbb{L}_i)\in H^2(\bar\M_{g,h};\mathbb Q).
$$

Let us define
$$
\Lambda_g^\vee(u)=u^g -\lambda_1 u^{g-1}+\cdots+(-1)^g \lambda_g.
$$
The Mari\~{n}o-Vafa conjecture states that the generating series
of the triple Hodge integrals
$$
\int_{\bar\M_{g,h}}\frac{\Lambda_g^\vee(1)\Lambda_g^\vee(\tau)
\Lambda_g^\vee(-\tau-1)}{\prod_{i=1}^h(1-\mu_i\psi_i)},
$$
for all $g$ and all $h$ can be expressed by close formulas of finite
expression in terms of the representations of symmetric groups, or
the Chern-Simons knot invariants. Here $\tau$ is a parameter and
$\mu_i$ are some integers. Many interesting Hodge integral
identities can be easily derived from this formula.

The Mari\~{n}o-Vafa conjecture originated from the large $N$ duality
between the Chern-Simons and string theory. It was proved by
exploring differential equations from both geometry and
combinatorics. The interested reader may read \cite{llz} for more
details.

\section{The Background of the Teichm\"uller Theory}

In this section, we recall some basic facts in Teichm\"uller theory
and introduce various notations for the following discussions.
Please see \cite{ga1} and \cite{tro1} for more details.

Let $\Sigma$ be an orientable surface with genus $g\geq 2$. A
complex structure on $\Sigma$ is a covering of $\Sigma$ by charts
such that the transition functions are holomorphic. By the
uniformization theorem, if we put a complex structure on $\Sigma$,
then it can be viewed as a quotient of the hyperbolic plane $\mathbb
H^2$ by a Fuchsian group. Thus there is a unique K\"ahler-Einstein
metric, or the hyperbolic metric on $\Sigma$.

Let $\mathcal C$ be the set of all complex structures on $\Sigma$.
Let $Diff^+(\Sigma)$ be the group of orientation preserving
diffeomorphisms and let $Diff^{+}_{0}(\Sigma)$ be the subgroup of
$Diff^+(\Sigma)$ consisting of those elements which are isotopic
to identity.

The groups $Diff^+(\Sigma)$ and $Diff^+_0(\Sigma)$ act naturally on
the space $\mathcal C$ by pull-back. The Teichm\"uller space is a
quotient of the space $\mathcal C$
\[
\mathcal T_g=\mathcal C /Diff^+_0(\Sigma).
\]
 From the famous Bers embedding theorem, now we know
that $\mathcal T_g$ can be embedded into $\mathbb C^{3g-3}$ as a
pseudoconvex domain and is contractible. Let
\[
\text{Mod}_g=Diff^+(\Sigma) /Diff^+_0(\Sigma)
\]
be the group of isotopic classes of diffeomorphisms. This group is
called the (Teichm\"uller) moduli group or the mapping class group.
 Its representations are of great interests in topology and in
 quantum field theory.

The moduli space $\M_g$ is the space of distinct complex structures on
$\Sigma$ and is defined to be
\[
\M_g=\mathcal C/Diff^+(\Sigma)=\mathcal T_g/\text{Mod}_g.
\]
The moduli space is a complex orbifold.

For any point $s\in \mathcal M_g$, let $X=X_s$ be a representative of the
corresponding class of Riemann surfaces. By the Kodaira-Spencer deformation
theory and the Hodge theory, we have
\[
T_X \mathcal M_g\cong H^1(X,T_X)=HB(X)
\]
where $HB(X)$ is the space of harmonic Beltrami differentials on
$X$.
\[
T_X^\ast \mathcal M_g\cong Q(X)
\]
where $Q(X)$ is the space of holomorphic quadratic differentials
on $X$.

Pick $\mu \in HB(X)$ and $\phi\in Q(X)$. If we fix a holomorphic
local coordinate $z$ on $X$, we can write
$\mu=\mu(z)\frac{\partial}{\partial z}\otimes d\bar z$ and
$\phi=\phi(z)dz^2$. Thus the duality between $T_X \mathcal M_g$
and $T_X^\ast \mathcal M_g$ is
\[
[\mu:\phi]=\int_X \mu(z)\phi(z)dzd\bar z.
\]

By the Riemann-Roch theorem, we have
\[
\dim_{\mathbb C}HB(X)=\dim_{\mathbb C}Q(X)=3g-3
\]
which implies
\[
\dim_{\mathbb C}\mathcal T_g=\dim_{\mathbb C}\M_g=3g-3.
\]

\section{Metrics on the Teichm\"uller Space and the Moduli Space}
There are many very famous classical metrics on the Teichm\"uller
and the moduli spaces and they have been studied independently by
many famous mathematicians. Each metric has played important role in
the study of the geometry and topology of the moduli and
Teichm\"uller spaces.

There are three Finsler metrics: the Teichm\"uller metric
$\Vert\cdot\Vert_T$, the Kobayashi metric $\Vert\cdot\Vert_K$ and
the Carath\'eodory metric $\Vert\cdot\Vert_C$. They are all
complete metrics on the Teichm\"uller space and are invariant
under the moduli group action. Thus they descend down to the
moduli space as complete Finsler metrics.

There are seven  K\"ahler metrics: the Weil-Petersson metric
$\omega_{_{WP}}$ which is incomplete, the Cheng-Yau's
K\"ahler-Einstein metric $\omega_{_{KE}}$, the McMullen metric
$\omega_{_{C}}$, the Bergman metric $\omega_{_{B}}$, the asymptotic
Poincar\'e metric on the moduli space $\omega_{_{P}}$, the Ricci
metric $\omega_{\tau}$ and the perturbed Ricci metric
$\omega_{\widetilde\tau}$. The last six metrics are complete. The
last two metrics are new metrics studied in details in \cite{lsy1}
and \cite{lsy2}.

Now let us give the precise definitions of these metrics and state
their basic properties.

The Teichm\"uller metric was first introduced by Teichm\"uller as
the $L^1$ norm in the cotangent space. For each
$\phi=\phi(z)dz^2\in Q(X)\cong T_X^\ast\M_g$, the Teichm\"uller
norm of $\phi$ is
\[
\Vert\phi\Vert_T=\int_X |\phi(z)|\ dzd\bar z.
\]
By using the duality, for each $\mu\in HB(X)\cong T_X\M_g$,
\[
\Vert\mu\Vert_T=\sup\{Re[\mu;\phi]\mid \Vert\phi\Vert_T=1\}.
\]
Please see \cite{ga1} for details. It is known that Teichm\"uller
metric has constant holomorphic sectional curvature $-1$.

The Kobayashi and the Carath\'eodory metrics can be defined for any
complex space in the following way: Let $Y$ be a complex manifold
and of dimension $n$. let $\Delta_R$ be the disk in $\mathbb C$ with
radius $R$. Let $\Delta=\Delta_1$ and let $\rho$ be the Poincar\'e
metric on $\Delta$. Let $p\in Y$ be a point and let $v\in T_p Y$ be
a holomorphic tangent vector. Let $\text{Hol}(Y,\Delta_R)$ and
$\text{Hol}(\Delta_R,Y)$ be the spaces of holomorphic maps from $Y$
to $\Delta_R$ and from $\Delta_R$ to $Y$ respectively. The
Carath\'eodory norm of the vector $v$ is defined to be
\[
\Vert v\Vert_C=\sup_{f\in\text{Hol}(Y,\Delta)}\Vert f_\ast
v\Vert_{\Delta,\rho}
\]
and the Kobayashi norm of $v$ is defined to be
\[
\Vert v\Vert_K=\inf_{f\in\text{Hol}(\Delta_R,Y),\ f(0)=p,\
f'(0)=v}\frac{2}{R}.
\]

The Bergman (pseudo) metric can also be defined for any complex
space $Y$ provided the Bergman kernel is positive. Let $K_Y$ be
the canonical bundle of $Y$ and let $W$ be the space of $L^2$
holomorphic sections of $K_Y$ in the sense that if $\sigma\in W$,
then
\[
\Vert\sigma\Vert_{L^2}^2=\int_Y
(\sqrt{-1})^{n^2}\sigma\wedge\bar\sigma<\infty.
\]
The inner product on $W$ is defined to be
\[
(\sigma,\rho)=\int_Y (\sqrt{-1})^{n^2}\sigma\wedge\bar\rho
\]
for all $\sigma,\rho\in W$. Let $\sigma_1,\sigma_2,\cdots$ be an
orthonormal basis of $W$. The Bergman kernel form is the
non-negative $(n,n)$-form
\[
B_Y=\sum_{j=1}^\infty(\sqrt{-1})^{n^2}\sigma_j\wedge\bar\sigma_j.
\]

With a choice of local coordinates $z_i,\cdots,z_n$, we have
\[
B_Y=BE_Y(z,\bar z)(\sqrt{-1})^{n^2}dz_1\wedge\cdots\wedge dz_n
\wedge d\bar z_1\wedge\cdots\wedge d\bar z_n
\]
where $BE_Y(z,\bar z)$ is called the Bergman kernel function. If
the Bergman kernel $B_Y$ is positive, one can define the Bergman
metric
\[
B_{i\bar j}=\frac{\partial^2\log BE_Y(z,\bar z)}{\partial z_i
\partial \bar z_j}.
\]
The Bergman metric is well-defined and is nondegenerate if the
elements in $W$ separate points and the first jet of $Y$. In this
case, the Bergman metric is a K\"ahler metric.

\begin{rem}
Both the Teichm\"uller space and the moduli space are equipped with
the Bergman metrics. However, the Bergman metric on the moduli space
is different from the metric induced from the Bergman metric of the
Teichm\"uller space. The Bergman metric defined on the moduli space
is incomplete due to the fact that the moduli space is
quasi-projective and any $L^2$ holomorphic section of the canonical
bundle can be extended over. However, the induced one is complete as
we shall see later.
\end{rem}

The basic properties of the Kobayashi, the Carath\'eodory and the
Bergman metrics are stated in the following proposition. Please see
\cite{ko2} for the details.
\begin{prop}
Let $X$ be a complex space. Then
\begin{enumerate}
\item $\Vert \cdot\Vert_{C,X}\leq \Vert \cdot \Vert_{K,X}$; \item
Let $Y$ be another complex space and $f:X\to Y$ be a holomorphic
map. Let $p\in X$ and $v\in T_p X$. Then $\Vert
f_\ast(v)\Vert_{C,Y,f(p)}\leq \Vert v \Vert_{C,X,p}$ and $\Vert
f_\ast(v)\Vert_{K,Y,f(p)}\leq \Vert v \Vert_{K,X,p}$; \item If
$X\subset Y$ is a connected open subset and $z\in X$ is a point.
Then with any local coordinates we have $BE_Y(z)\leq BE_X(z)$;
\item If the Bergman kernel is positive, then at each point $z\in
X$, a peak section $\sigma$ at $z$ exists. Such a peak section is
unique up to a constant factor $c$ with norm $1$. Furthermore,
with any choice of local coordinates, we have
$BE_X(z)=|\sigma(z)|^2$; \item If the Bergman kernel of $X$ is
positive, then $\Vert \cdot \Vert_{C,X}\leq 2\Vert \cdot
\Vert_{B,X}$; \item If $X$ is a bounded convex domain in $\mathbb
C^n$, then $\Vert \cdot\Vert_{C,X}= \Vert \cdot \Vert_{K,X}$;
\item Let $|\cdot |$ be the Euclidean norm and let $B_r$ be the
open ball with center $0$ and radius $r$ in $\mathbb C^n$. Then
for any holomorphic tangent vector $v$ at $0$,
\[
\Vert v\Vert_{C,B_r,0}=\Vert v\Vert_{K,B_r,0}=\frac{2}{r}|v|
\]
where $|v|$ is the Euclidean norm of $v$.
\end{enumerate}
\end{prop}

The three Finsler metrics have been very powerful tools in
understanding the hyperbolic geometry of the moduli spaces, and the
mapping class group. It is also known since 70's that the Bergman
metric on the Teichm\"uller space is complete.

The Weil-Petersson metric is the first K\"ahler metric defined on
the Teichm\"uller and the moduli space. It is defined by using the
$L^2$ inner
product on the tangent space in the following way:\\
Let $\mu,\nu\in T_X\M_g$ be two tangent vectors and let $\lambda$
be the unique K\"ahler-Einstein metric on $X$. Then the
Weil-Pertersson metric is
\[
h(\mu,\nu)=\int_X \mu\bar\nu\ dv
\]
where $dv=\frac{\sqrt{-1}}{2}\lambda dz\wedge d\bar z$ is the
volume form. Details can be found in \cite{lsy1}, \cite{ma1} and
\cite{wol1}.

The curvatures of the Weil-Petersson metric have been
well-understood due to the works of Ahlfors, Royden and Wolpert. Its
Ricci and holomorphic sectional curvature are all negative with
negative upper bound, but with no lower bound. Its boundary behavior
is understood, from which it is not hard to see that it is an
incomplete metric.

The existence of the K\"ahler-Einstein metric was given by the
work of Cheng-Yau \cite{cy2}. Its Ricci curvature is
$-1$. Namely,
\[
\partial\bar\partial\log\omega_{_{KE}}^{n}=\omega_{_{KE}}
\]
where $n=3g-3$. They actually proved that a bounded domain in
$\mathbb{C}^n$ admits a complete K\"ahler-Einstein metric if and
only if it is pseudoconvex.

The McMullen $1/l$ metric defined in \cite{mc} is a perturbation of
the Weil-Petersson metric by adding a K\"ahler form whose potential
involves the short geodesic length functions on the Riemann
surfaces. For each simple closed curve $\gamma$ in $X$, let
$l_\gamma(X)$ be the length of the unique geodesic in the homotopy
class of $\gamma$ with respect to the unique K\"ahler-Einstein
metric. Then the McMullen metric is defined as
\[
\omega_{1/l}=\omega_{_{WP}}-i\delta\sum_{l_\gamma(X)<\epsilon}
\partial\bar\partial\text{Log}\frac{\epsilon}{l_\gamma}
\]
where $\epsilon$ and $\delta$ are small positive constants and
$\text{Log}(x)$ is a smooth function defined as
\[
\text{Log}(x)= \begin{cases} \log x & x\geq 2\\
0 & x\leq 1.
\end{cases}
\]
This metric is K\"ahler hyperbolic which means it satisfies the
following conditions:
\begin{enumerate}
\item $(\M_g,\omega_{1/l})$ has finite volume;
\item The sectional curvature of $(\M_g,\omega_{1/l})$ is bounded
above and below;
\item The injectivity radius of  $(\mathcal T_g,\omega_{1/l})$ is
bounded below;
\item On $\mathcal T_g$, the K\"ahler form $\omega_{1/l}$ can be
written as $\omega_{1/l}=d\alpha$ where $\alpha$ is a bounded
$1$-form.
\end{enumerate}
An immediate consequence of the K\"ahler hyperbolicity is that the
$L^2$-cohomology is trivial except for the middle dimension.

The asymptotic Poincar\'e metric can be defined as a complete
K\"ahler metric on a complex manifold $M$ which is obtained by
removing a divisor $Y$ with only normal crossings from a compact
K\"ahler manifold $(\bar M,\omega)$.

Let $\bar M$ be a compact K\"ahler manifold of dimension $m$. Let
$Y\subset \bar M$ be a divisor of normal crossings and let $M=\bar
M\setminus Y$. Cover $\bar M$ by coordinate charts
$U_1,\cdots,U_p,\cdots,U_q$ such that $(\bar
U_{p+1}\cup\cdots\cup\bar U_q)\cap Y=\Phi$. We also assume that
for each $1\leq \alpha \leq p$, there is a constant $n_\alpha$
such that $U_\alpha\setminus
Y=(\Delta^\ast)^{n_\alpha}\times\Delta^{m-n_\alpha}$ and on
$U_\alpha$, $Y$ is given by $z_1^\alpha\cdots
z_{n_\alpha}^\alpha=0$. Here $\Delta$ is the disk of radius
$\frac{1}{2}$ and $\Delta^\ast$ is the punctured disk of radius
$\frac{1}{2}$. Let $\{\eta_i\}_{1\leq i\leq q}$ be the partition
of unity subordinate to the cover $\{U_i\}_{1\leq i\leq q}$. Let
$\omega$ be a K\"ahler metric on $\bar M$ and let $C$ be a
positive constant. Then for $C$ large, the K\"ahler form
\[
\omega_p=C\omega+\sum_{i=1}^p\sqrt{-1}\partial\bar\partial \bigg
(\eta_i\log\log\frac{1}{z_1^i\cdots z_{n_i}^i}\bigg )
\]
defines a complete metric on $M$ with finite volume since on each
$U_i$ with $1\leq i\leq p$, $\omega_p$ is bounded from above and
below by the local Poincar\'e metric on $U_i$. We call this metric
the asymptotic Poincar\'e metric.

The signs of the curvatures of the above metrics are all unknown. We
actually only know that the K\"ahler-Einstein metric has constant
negative Ricci curvature and that the McMullen metric has bounded
geometry. Also except the asymptotic Poincar\'e metric, the boundary
behaviors of the other metrics are unknown either before our works.
It is interesting that to understand them we need to introduce new
metrics.

Now we define the Ricci metric and the perturbed Ricci metric. The
curvature properties and asymptotics of these two new metrics are
understood by us and will be stated in the following sections.
Please also see \cite{lsy1} and \cite{lsy2} for details.

By the works of Ahlfors, Royden and Wolpert we know that the Ricci
curvature of the Weil-Petersson metric has a negative upper bound.
Thus we can use the negative Ricci form of the Weil-Petersson metric
as the K\"ahler form of a new metric. We call this metric the Ricci
metric and denote it by $\tau$. That is
\[
\omega_\tau=\partial\bar\partial\log\omega_{_{WP}}^n.
\]
Through careful analysis, we now understand that the Ricci metric is
a natural canonical complete K\"ahler metric with many good
properties. However, its holomorphic sectional curvature is only
asymptotically negative. To get a metric with good sign on its
curvatures, we introduced the perturbed Ricci metric
$\omega_{\widetilde\tau}$ as a combination of the Ricci metric and
the Weil-Petersson metric:
\[
\omega_{\widetilde\tau}=\omega_\tau+C\omega_{_{WP}}
\]
where $C$ is a large positive constant. As we shall see later that
the perturbed Ricci metric has desired curvature properties so that
we can put it either on the target or on the domain manifold in
Yau's Schwarz lemma, from which we can compare the above metrics.

\section{The Curvature Formulas}

In this section we describe the harmonic lift of a vector field on
the moduli space to the universal curve due to Royden, Siu
\cite{siu1} and Schumacher \cite{sc1}. Details can also be found in
\cite{lsy1}. We then use this method to derive the curvature formula
for the Weil-Petersson metric, the Ricci metric and the perturbed
Ricci metric.

To compute the curvature of a metric on the moduli space, we need to
take derivatives of the metric in the direction of the moduli space.
However, it is quite difficult to estimate the curvature by using a
formula obtained in such a way. The central idea is to obtain a
formula where the derivatives are taken in the fiber direction. We
can view the deformation of complex structures on a topological
surface as the deformation of the K\"ahler-Einstein metrics.

Let $\M_g$ be the moduli space of Riemann surfaces of genus $g$
where $g \geq 2$. Let $n=3g-3$ be the complex dimension of
$\M_g$. Let $\frak X$ be the total space and let $\pi:\frak X \to
\M_g$ be the projection map.

Let $s_1, \cdots, s_n$ be holomorphic local coordinates near a
regular point $s \in \M_g$ and assume that $z$ is a holomorphic
local  coordinate on the fiber $X_s=\pi^{-1}(s)$. For holomorphic
vector fields $\frac{\pa}{\pa s_1},\cdots, \frac{\pa}{\pa s_n}$,
there are vector fields $v_1,\cdots, v_n$ on $\frak X$ such that
\begin{enumerate}
\item  $\pi_*(v_i)=\frac{\pa}{\pa s_i}$ for $i=1,\cdots,n$; \item
$\bar\pa v_i$ are harmonic $TX_s$-valued $(0,1)$ forms for
$i=1,\cdots,n$.
\end{enumerate}
The vector fields $v_1,\cdots,v_n$ are called the harmonic lift of
the vectors $\frac{\pa}{\pa s_1},\cdots, \frac{\pa}{\pa s_n}$. The
existence of such harmonic vector fields was pointed by Siu.
Schumacher in his work gave an explicit construction of such lift.
We now describe it.

Since $g \geq 2$, we can assume that each fiber is equipped with
the K\"ahler-Einstein metric
$\lambda=\frac{\sqrt{-1}}{2}\lambda(z,s)dz\wedge d\bar z$. The
K\"ahler-Einstein condition gives the following equation:
\begin{eqnarray}\label{100}
\pz\pzb\log\lambda=\lambda.
\end{eqnarray}
For the rest of this paper we denote $\frac{\partial}{\partial
s_i}$ by $\partial_i$  and $\frac{\pa}{\pa z}$ by $\pa_z$. Let
\begin{eqnarray*}
a_i=-\lambda^{-1}\partial_i\pzb\log\lambda
\end{eqnarray*}
and let
\begin{eqnarray*}
A_i=\pzb a_i.
\end{eqnarray*}
Then the harmonic horizontal lift of $\partial_i$ is
\[
v_i=\partial_i+a_i\pz.
\]
In particular
\[
B_i=A_i\pz\otimes d\bar z \in H^{1}(X_s,T_{X_s})
\]
is harmonic. Furthermore, the lift $\pa_i\mapsto B_i$ gives the
Kodaira-Spencer map $T_s\M_g\rightarrow H^1(X_s,T_{X_s})$. Thus the
Weil-Petersson metric on $\M_g$ is
\begin{eqnarray*}
h_{i\bar j}(s)=\int_{X_s}B_i \cdot \bar{B_j} \ dv=
\int_{X_s}A_i\bar{A_j}\ dv,
\end{eqnarray*}
where $dv=\frac{\sqrt{-1}}{2}\lambda dz\wedge d\bar z$ is the
volume form on the fiber $X_s$.

Let $R_{i\bar jk\bar l}$ be the curvature tensor of the
Weil-Petersson metric. Here we adopt the following notation for
the curvature of a K\"ahler metric:\\
For a K\"ahler metric $(M,g)$, the curvature tensor is given by
\[
R_{i\bar j k\bar l}=\frac{\partial^2 g_{i\bar j}}{\partial z_k
\partial \bar z_l}-g^{p\bar q}\frac{\partial g_{i\bar q}}
{\partial z_k}\frac{\partial g_{p\bar j}} {\partial \bar z_l}.
\]
In this case, the Ricci curvature is given by
\[
R_{i\bar j}=-g^{k\bar l}R_{i\bar j k\bar l}.
\]

By using the curvature of the Weil-Petersson metric, we can define
the Ricci metric:
\[
\tau_{i\bar j}=h^{k\bar l}R_{i\bar j k\bar l}
\]
and the perturbed Ricci metric:
\[
\widetilde\tau_{i\bar j}=\tau_{i\bar j}+Ch_{i\bar j}
\]
where $C$ is a positive constant.

Before we present the curvature formulas for the above metrics, we
need to introduce the Maass operators and norms on a Riemann surface
\cite{wol1}.

Let $X$ be a Riemann surface and let $\kappa$ be its canonical
bundle. For any integer $p$, let $S(p)$ be the space of smooth
sections of $(\kappa\otimes\bar\kappa^{-1})^{\frac{p}{2}}$. Fix a
conformal metric $ds^2=\rho^2(z)|dz|^2$. In the following, we will
take $ds^2$ to be the K\"ahler-Einstein metric although the
following definitions work for all metrics.

The Maass operators $K_p$ and $L_p$ are defined to be the metric
derivatives \\
$K_p:S(p)\to S(p+1)$ and $L_p:S(p)\to S(p-1)$ given by
\[
K_p(\sigma)=\rho^{p-1}\pz(\rho^{-p}\sigma)
\]
and
\[
L_p(\sigma)=\rho^{-p-1}\pzb(\rho^{p}\sigma)
\]
where $\sigma \in S(p)$.

The operators $P=K_1K_0$ and $\Box=-L_1K_0$ will play important
roles in the curvature formulas. Here the operator $\Box$ is just
the Laplace operator. We also let $T=(\Box+1)^{-1}$ to be the
Green operator.

Each element $\sigma \in S(p)$ have a well-defined absolute value
$|\sigma|$ which is independent of the choice of local coordinate.
We define the $C^k$ norm of $\sigma$:\\
Let $Q$ be an operator
which is a composition of operators $K_\ast$ and $L_\ast$. Denote
by $|Q|$ the number of factors. For any $\sigma \in S(p)$, define
\[
\Vert\sigma\Vert_{0}=\sup_{X}|\sigma|
\]
and
\[
\Vert\sigma\Vert_k=\sum_{|Q|\leq k}\Vert Q\sigma\Vert_0.
\]
We can also localize the norm on a subset of $X$. Let $\Omega
\subset X$ be a domain. We can define
\[
\Vert\sigma\Vert_{0,\Omega}=\sup_{\Omega}|\sigma|
\]
and
\[
\Vert\sigma\Vert_{k,\Omega} =\sum_{|Q|\leq k}\Vert
Q\sigma\Vert_{0,\Omega}.
\]

We let $f_{i\bar j}=A_i\bar A_j$ and $e_{i\bar j}=T(f_{i\bar j})$.
These functions will be the building blocks for the curvature
formulas.

The trick of converting derivatives from the moduli directions to
the fiber directions is the following lemma due to Siu and
Schumacher:
\begin{lemma}
Let $\eta$ be a relative $(1,1)$-form on the total space $\frak
X$. Then
\[
\frac{\partial}{\partial s_i}\int_{X_s}\eta=\int_{X_s}L_{v_i}\eta.
\]
\end{lemma}

The curvature formula of the Weil-Petersson metric was first
established by Wolpert by using a different method \cite{wol3} and
later was generalized by Siu \cite{siu1} and Schumacher \cite{sc1}
by using the above lemma:
\begin{theorem}
The curvature of the Weil-Petersson metric is given by
\begin{eqnarray}\label{200}
R_{i\bar j k\bar l}=\int_{X_s}(e_{i\bar j}f_{k\bar l}+ e_{i\bar
l}f_{k\bar j})\ dv.
\end{eqnarray}
\end{theorem}

For the proof, please see \cite{lsy1}. From this formula it is
rather easy to show that the Ricci and the holomorphic sectional
curvature have explicit negative upper bound.

To establish
the curvature formula of the Ricci metric, we need to introduce
more operators.
Firstly, the commutator of the operator $v_k$ and $(\Box+1)$ will
play an important role. Here we view the vector field $v_k$ as a
operator acting on functions. We define
\[
\xi_k =[\Box+1,v_k].
\]
A direct computation shows that
\[
\xi_k=-A_k P.
\]

Also we can define the commutator of $\bar{v_l}$ and $\xi_k$. Let
\[
Q_{k\bar l}=[\bar{v_l}, \xi_k].
\]
We have
\[
Q_{k\bar l}(f)=\bar P(e_{k\bar l})P(f)-2f_{k\bar l}\Box f
+\lambda^{-1}\pz f_{k\bar l}\pzb f
\]
for any smooth function $f$.

To simplify the notation, we introduce the symmetrization operator
of the indices. Let $U$ be any quantity which depends on indices
$i,k,\alpha,\bar j,\bar l, \bar\beta$. The symmetrization operator
$\sigma_1$ is defined by taking summation of all orders of the triple
$(i,k,\alpha)$. That is
\begin{align*}
\begin{split}
\sigma_1(U(i,k,\alpha,\bar j,\bar l, \bar\beta))=&
U(i,k,\alpha,\bar j,\bar l, \bar\beta)+ U(i,\alpha,k,\bar j,\bar
l, \bar\beta)+ U(k,i,\alpha,\bar j,\bar l, \bar\beta)+
U(k,\alpha,i,\bar j,\bar l, \bar\beta)\\
&+U(\alpha,i,k,\bar j,\bar l, \bar\beta)+ U(\alpha,k,i,\bar j,\bar
l, \bar\beta).
\end{split}
\end{align*}
Similarly, $\sigma_2$ is the symmetrization operator of $\bar j$
and $\bar \beta$ and $\widetilde{\sigma_1}$ is the symmetrization
operator of $\bar j$, $\bar l$ and $\bar \beta$.

In \cite{lsy1} the following curvature formulas for the Ricci and
perturbed Ricci metric were proved:
\begin{theorem}\label{riccicurv}
Let $s_1,\cdots,s_n$ be local holomorphic coordinates at $s \in
M_g$. Then at $s$, we have
\begin{align}\label{finalcurv}
\begin{split}
\widetilde{R}_{i\bar j k\bar l}=&h^{\alpha\bar\beta}
\left\{\sigma_1\sigma_2\int_{X_s}
\left\{(\Box+1)^{-1}(\xi_k(e_{i\bar j}))
\bar{\xi}_l(e_{\alpha\bar\beta})+ (\Box+1)^{-1} (\xi_k(e_{i\bar
j})) \bar{\xi}_\beta(e_{\alpha\bar l})
\right\}\ dv\right\}\\
&+h^{\alpha\bar\beta} \left\{\sigma_1\int_{X_s}Q_{k\bar
l}(e_{i\bar j}) e_{\alpha\bar\beta}\ dv
\right\}\\
&-\tau^{p\bar q}h^{\alpha\bar\beta}h^{\gamma\bar\delta}
\left\{\sigma_1\int_{X_s}\xi_k(e_{i\bar q}) e_{\alpha\bar\beta}\
dv\right\}\left\{ \widetilde\sigma_1\int_{X_s}\bar{\xi}_l(e_{p\bar
j})
e_{\gamma\bar\delta})\ dv\right\}\\
&+\tau_{p\bar j}h^{p\bar q}R_{i\bar q k\bar l}
\end{split}
\end{align}
and
\begin{align}\label{finalpercurv}
\begin{split}
P_{i\bar j k\bar l}=&h^{\alpha\bar\beta}
\left\{\sigma_1\sigma_2\int_{X_s}
\left\{(\Box+1)^{-1}(\xi_k(e_{i\bar j}))
\bar{\xi}_l(e_{\alpha\bar\beta})+ (\Box+1)^{-1} (\xi_k(e_{i\bar
j})) \bar{\xi}_\beta(e_{\alpha\bar l})
\right\}\ dv\right\}\\
&+h^{\alpha\bar\beta} \left\{\sigma_1\int_{X_s}Q_{k\bar
l}(e_{i\bar j}) e_{\alpha\bar\beta}\ dv
\right\}\\
&-\widetilde\tau^{p\bar q}h^{\alpha\bar\beta}h^{\gamma\bar\delta}
\left\{\sigma_1\int_{X_s}\xi_k(e_{i\bar q}) e_{\alpha\bar\beta}\
dv\right\}\left\{ \widetilde\sigma_1\int_{X_s}\bar{\xi}_l(e_{p\bar
j})
e_{\gamma\bar\delta})\ dv\right\}\\
&+\tau_{p\bar j}h^{p\bar q}R_{i\bar q k\bar l} +CR_{i\bar j k\bar
l}.
\end{split}
\end{align}
where $R_{i\bar j k\bar l}$, $\widetilde R_{i\bar j k\bar l}$, and
$P_{i\bar j k\bar l}$ are the curvature of the Weil-Petersson
metric, the Ricci metric and the perturbed Ricci metric
respectively.
\end{theorem}

Unlike the curvature formula of the Weil-Petersson metric which we
can see the sign of the curvature directly, the above formulas are
too complicated and we cannot see the sign. So we need to study the
asymptotic behaviors of these curvatures, and first the metrics
themselves.

\section{The Asymptotics}

To compute the asymptotics of these metrics and their curvatures, we
first need to find a canonical way to construct local coordinates
near the boundary of the moduli space. We first describe the
Deligne-Mumford compactification of the moduli space and introduce
the pinching coordinate and the plumbing construction which due to
Earle and Marden. Please see \cite{ma1}, \cite{wol1}, \cite{tr1} and
\cite{lsy1} for details.

A point $p$ in a Riemann surface $X$ is a node if there
is a neighborhood of $p$ which is isometric to the
germ $\{ (u,v)\mid uv=0,\ |u|,|v|<1 \} \subset
\mathbb{C}^2$. Let $p_1,\cdots,p_k$ be the nodes on $X$. $X$ is called
stable if each connected component of $X\setminus \{p_1,\cdots,p_k\}$
has negative Euler characteristic. Namely, each connected component
has a unique complete hyperbolic metric.

Let $\mathcal M_g$ be the moduli space of Riemann
surfaces of genus $g \geq 2$. The Deligne-Mumford compactification
$\bar\M_g$ is the union of $\mathcal M_g$ and corresponding stable
nodal surfaces \cite{dm1}. Each point $y \in \bar{\mathcal M}_g
\setminus \mathcal M_g$ corresponds to a stable noded
surface $X_y$.

We recall the rs-coordinate on a
Riemann surface defined by Wolpert in \cite{wol1}.
There are two cases: the puncture case and the short
geodesic case. For the puncture case, we have a noded
surface $X$ and a node $p\in X$. Let $a,b$ be two punctures
which are paired to form $p$.
\begin{definition}
The local coordinate charts $(U,u)$ near $a$ is called rs-coordinate
if $u(a)=0$, $u$ maps $U$ to the punctured disc $0<|u|<c$ with $c>0$
and the restriction to $U$ of the K\"ahler-Einstein metric on $X$
can be written as $\frac{1}{2|u|^2(\log |u|)^2} |du|^2$. The
rs-coordinate $(V,v)$ near $b$ is defined in a similar way.
\end{definition}
For the short geodesic case, we have a closed surface
$X$, a closed geodesic $\gamma \subset X$ with length
$l <c_\ast$ where $c_\ast$ is the collar constant.
\begin{definition}
The local coordinate chart $(U,z)$ is called rs-coordinate at
$\gamma$ if $\gamma \subset U$, $z$ maps $U$ to the annulus
$c^{-1}|t|^{\frac{1}{2}}<|z| <c|t|^{\frac{1}{2}}$ and the
K\"ahler-Einstein metric on $X$ can be written as
$\frac{1}{2}(\frac{\pi}{\log |t|}\frac{1}{|z|}\csc \frac{\pi\log
|z|}{\log |t|})^2 |dz|^2$.
\end{definition}
\begin{rem}
We put the factor $\frac{1}{2}$ in the above two
definitions to normalize such that \eqref{100} holds.
\end{rem}

By Keen's collar theorem \cite{ke1}, we have the
following lemma:
\begin{lemma}\label{gcollar}
Let $X$ be a closed surface and let $\gamma$ be a
closed geodesic on $X$ such that the length $l$
of $\gamma$ satisfies $l <c_\ast$. Then there is
a collar $\Omega$ on $X$ with holomorphic coordinate
$z$ defined on $\Omega$ such that
\begin{enumerate}
\item $z$ maps $\Omega$ to the annulus
$\frac{1}{c}e^{-\frac{2\pi^2}{l}}<|z|<c$ for $c>0$;
\item the K\"ahler-Einstein metric on $X$ restrict to
$\Omega$ is given by
\begin{eqnarray}\label{precmetric}
(\frac{1}{2}u^2 r^{-2}\csc^2\tau) |dz|^2
\end{eqnarray}
where $u=\frac{l}{2\pi}$, $r=|z|$ and $\tau=u\log r$;
\item the geodesic $\gamma$ is given by $|z|=
e^{-\frac{\pi^2}{l}}$.
\end{enumerate}
We call such a collar $\Omega$ a genuine collar.
\end{lemma}
We notice that the constant $c$ in the above lemma has
a lower bound such that the area of $\Omega$ is bounded
from below. Also, the coordinate $z$ in the above
lemma is rs-coordinate. In the following, we will keep
the notation $u$, $r$ and $\tau$.

Now we describe the local manifold cover of $\bar{\mathcal M}_g$
near the boundary. We take the construction of Wolpert \cite{wol1}.
Let $X_{0,0}$ be a noded surface corresponding to a codimension $m$
boundary point. $X_{0,0}$ have $m$ nodes $p_1,\cdots,p_m$.
$X_0=X_{0,0}\setminus \{ p_1,\cdots,p_m \}$ is a union of punctured
Riemann surfaces. Fix rs-coordinate charts $(U_i,\eta_i)$ and
$(V_i,\zeta_i)$ at $p_i$ for $i=1,\cdots,m$ such that all the $U_i$
and $V_i$ are mutually disjoint. Now pick an open set $U_0 \subset
X_0$ such that the intersection of each connected component of $X_0$
and $U_0$ is a nonempty relatively compact set and the intersection
$U_0 \cap (U_i\cup V_i)$ is empty for all $i$. Now pick Beltrami
differentials $\nu_{m+1},\cdots,\nu_{n}$ which are supported in
$U_0$ and span the tangent space at $X_0$ of the deformation space
of $X_0$. For $s=(s_{m+1},\cdots,s_n)$, let $\nu(s)=\sum_{i=m+1}^n
s_i\nu_i$. We assume $|s|=(\sum |s_i|^2)^{\frac{1}{2}}$ is small
enough such that $|\nu(s)|<1$. The noded surface $X_{0,s}$ is
obtained by solving the Beltrami equation $\bar\partial
w=\nu(s)\partial w$. Since $\nu(s)$ is supported in $U_0$,
$(U_i,\eta_i)$ and $(V_i,\zeta_i)$ are still holomorphic coordinates
on $X_{0,s}$. Note that they are no longer rs-coordinates. By the
theory of Alhfors and Bers \cite{ab1} and Wolpert \cite{wol1} we can
assume that there are constants $\delta,c>0$ such that when
$|s|<\delta$, $\eta_i$ and $\zeta_i$ are holomorphic coordinates on
$X_{0,s}$ with $0<|\eta_i|<c$ and $0<|\zeta_i|<c$. Now we assume
$t=(t_1,\cdots,t_m)$ has small norm. We do the plumbing construction
on $X_{0,s}$ to obtain $X_{t,s}$. Remove from $X_{0,s}$ the discs
$0<|\eta_i|\leq \frac{|t_i|}{c}$ and $0<|\zeta_i|\leq
\frac{|t_i|}{c}$ for each $i=1,\cdots,m$. Now identify
$\frac{|t_i|}{c}<|\eta_i|< c$ with $\frac{|t_i|}{c}<|\zeta_i|< c$ by
the rule $\eta_i \zeta_i=t_i$. This defines the surface $X_{t,s}$.
The tuple $(t_1,\cdots,t_m,s_{m+1},\cdots,s_n)$ are the local
pinching coordinates for the manifold cover of $\bar{\mathcal M}_g$.
We call the coordinates $\eta_i$ (or $\zeta_i$) the plumbing
coordinates on $X_{t,s}$ and the collar defined by
$\frac{|t_i|}{c}<|\eta_i|< c$ the plumbing collar.
\begin{rem}
By the estimate of Wolpert \cite{wol2}, \cite{wol1} on
the length of short geodesic, the quantity
$u_i=\frac{l_i}{2\pi}\sim -\frac{\pi}{\log|t_i|}$.
\end{rem}

Now we describe the estimates of the asymptotics of these metrics
and their curvatures. The principle is that, when we work on a
nearly degenerated surface, the geometry focuses on the collars. Our
curvature formulas depend on the K\"ahler-Einstein metrics of the
family of Riemann surfaces near a boundary points. One can obtain
approximate K\"ahler-Einstein metric on these collars by the graft
construction of Wolpert \cite{wol1} which is done by gluing the
hyperbolic metric on the nodal surface with the model metric
described above.

To use the curvature formulas \eqref{200}, \eqref{finalcurv} and
\eqref{finalpercurv} to estimate the asymptotic behavior, one also
needs to analyze the transition from the plumbing coordinates on the
collars to the rs-coordinates. The harmonic Beltrami differentials
were constructed by Masur \cite{ma1} by using the plumbing
coordinates and it is easier to compute the integration by using
rs-coordinates. Such computation was done in \cite{tr1} by using the
graft metric of Wolpert and the maximum principle. A clear
description can be found in \cite{lsy1}. We have the following
theorem:
\begin{theorem}\label{imp}
Let $(t,s)$ be the pinching coordinates on $\bar{\mathcal M}_g$
near $X_{0,0}$ which corresponds to a codimension $m$ boundary point of
$\bar{\mathcal M}_g$. Then there exist constants $M,\delta>0$ and
$1>c>0$ such that if $|(t,s)|<\delta$, then the $j$-th plumbing
collar on $X_{t,s}$ contains the genuine collar $\Omega^j_c$.
Furthermore, one can choose rs-coordinate $z_j$ on the collar
$\Omega_c^j$ properly such that the holomorphic quadratic differentials
$\psi_1,\cdots,\psi_n$ corresponding to the cotangent vectors
$dt_1,\cdots,ds_n$ have form $\psi_i=\phi_i(z_j)dz_j^2$ on the
genuine collar $\Omega^j_c$ for $1\leq j \leq m$ where
\begin{enumerate}
\item $\phi_i(z_j)=\frac{1}{z_j^2}(q_i^j(z_j)+\beta_i^j)$ if $i\geq m+1$;
\item
  $\phi_i(z_j)=(-\frac{t_j}{\pi})\frac{1}{z_j^2}(q_j(z_j)+\beta_j)$ if
$i=j$;
\item $\phi_i(z_j)=(-\frac{t_i}{\pi})
\frac{1}{z_j^2}(q_i^j(z_j)+\beta_i^j)$ if $1\leq i \leq m$ and $i\ne j$.
\end{enumerate}
Here $\beta_i^j$ and $\beta_j$ are functions of $(t,s)$,
$q_i^j$ and $q_j$ are functions of $(t,s,z_j)$ given by
\[
q_i^j(z_j)=\sum_{k<0}\alpha_{ik}^j(t,s)t_j^{-k}z_j^k
+\sum_{k>0}\alpha_{ik}^j(t,s)z_j^k
\]
and
\[
q_j(z_j)=\sum_{k<0}\alpha_{jk}(t,s)t_j^{-k}z_j^k
+\sum_{k>0}\alpha_{jk}(t,s)z_j^k
\]
such that
\begin{enumerate}
\item $\sum_{k<0}|\alpha_{ik}^j|c^{-k}\leq M$ and
$\sum_{k>0}|\alpha_{ik}^j|c^{k}\leq M$ if $i\ne j$;
\item $\sum_{k<0}|\alpha_{jk}|c^{-k}\leq M$ and
$\sum_{k>0}|\alpha_{jk}|c^{k}\leq M$;
\item $|\beta_i^j|=O(|t_j|^{\frac{1}{2}-\epsilon})$
with $\epsilon<\frac{1}{2}$ if $i\ne j$;
\item $|\beta_j|=(1+O(u_0))$
\end{enumerate}
where $u_0=\sum_{i=1}^m u_i+\sum_{j=m+1}^n |s_j|$.
\end{theorem}

By definition, the metric on the cotangent bundle induced by the
Weil-Petersson metric is given by
\[
h^{i\bar j}=\int_{X_{t,s}}\lambda^{-2}\phi_i\bar\phi_j\ dv.
\]
We then have the following series of estimates, see \cite{lsy1}.
First by using this formula and taking inverse, we can estimate
the Weil-Petersson metric.
\begin{theorem}\label{wpasymp}
Let $(t,s)$ be the pinching coordinates. Then
\begin{enumerate}
\item $h^{i\bar i}=2u_i^{-3}|t_i|^2(1+O(u_0))$ and $h_{i\bar i}
=\frac{1}{2}\frac{u_i^{3}}{|t_i|^2}(1+O(u_0))$
for $1\leq i\leq m$;
\item $h^{i\bar j}=O(|t_it_j|)$ and $h_{i\bar
j}=O(\frac{u_i^3u_j^3}{|t_it_j|})$ if $1\leq i,j \leq m$ and $i\ne
j$;
\item $h^{i\bar j}=O(1)$ and $h_{i\bar j}=O(1)$ if $m+1\leq
i,j \leq n$;
\item $h^{i\bar j}=O(|t_i|)$ and $h_{i\bar
j}=O(\frac{u_i^3}{|t_i|})$ if $i\leq m < j$ or $j \leq m<i$ .
\end{enumerate}
\end{theorem}

Then we use the duality to construct the harmonic Beltrami
differentials. We have
\begin{lemma}\label{aj10}
On the genuine collar $\Omega_c^j$ for $c$ small, the coefficient
functions $A_i$ of the harmonic Beltrami differentials have the form:
\begin{enumerate}
\item
$A_i=\frac{z_j}{\bar{z_j}}\sin^2\tau_j
\big ( \bar{p_i^j(z_j)}+\bar{b_i^j}\big )$ if $i\ne j$;
\item
$A_j=\frac{z_j}{\bar{z_j}}\sin^2\tau_j(\bar{p_j(z_j)}+\bar{b_j})$
\end{enumerate}
where
\begin{enumerate}
\item $p_i^j(z_j)=\sum_{k\leq -1}a_{ik}^j\rho_j^{-k}z_j^k
+\sum_{k\geq 1}a_{ik}^jz_j^k$ if $i\ne j$;
\item  $p_j(z_j)=\sum_{k\leq -1}a_{jk}\rho_j^{-k}z_j^k
+\sum_{k\geq 1}a_{jk}z_j^k$.
\end{enumerate}
In the above expressions, $\rho_j=e^{-\frac{2\pi^2}{l_j}}$ and the
coefficients satisfy the following conditions:
\begin{enumerate}
\item $\sum_{k\leq -1}|a_{ik}^j|c^{-k}=O(u_j^{-2})$ and
$\sum_{k\geq 1}|a_{ik}^j|c^{k}=O(u_j^{-2})$ if $i\geq m+1$;
\item $\sum_{k\leq -1}|a_{ik}^j|c^{-k}=O(u_j^{-2})
O\big ( \frac{u_i^3}{|t_i|}\big )$ and
$\sum_{k\geq 1}|a_{ik}^j|c^{k}=O(u_j^{-2})
O\big ( \frac{u_i^3}{|t_i|}\big )$ if $i\leq m$ and $i\ne j$;
\item $\sum_{k\leq -1}|a_{jk}|c^{-k}=O(\frac{u_j}{|t_j|})$ and
$\sum_{k\geq 1}|a_{jk}|c^{k}=O(\frac{u_j}{|t_j|})$;
\item $|b_i^j|=O(u_j)$ if $i\geq m+1$;
\item $|b_i^j|=O(u_j)O\big ( \frac{u_i^3}{|t_i|}\big )$
if $i\leq m$ and $i\ne j$;
\item $b_j=-\frac{u_j}{\pi\bar{t_j}}(1+O(u_0))$.
\end{enumerate}
\end{lemma}

To use the curvature formulas to estimate the Ricci metric and the
perturbed Ricci metric, one needs to find accurate estimate of the
operator $T=(\Box+1)^{-1}$. More precisely, one needs to estimate
the functions $e_{i\bar j}=T(f_{i\bar j})$. To avoid writing down
the Green function of $T$, we construct approximate solutions and
localize on the collars in \cite{lsy1}. Pick a positive constant
$c_1<c$ and define the cut-off function $\eta\in C^{\infty}(\mathbb
R, [0,1])$ by
\begin{eqnarray}\label{cutoff}
\begin{cases}
\eta(x)=1 & x\leq \log c_1\\
\eta(x)=0 & x\geq \log c\\
0<\eta(x)<1 & \log c_1<x<\log c.
\end{cases}
\end{eqnarray}
It is clear that the derivatives of $\eta$ are bounded by constants which
only depend on $c$ and $c_1$. Let $\widetilde{e_{i\bar j}}(z)$ be the
function on $X$ defined in the following way where $z$ is taken to be
$z_i$ on the collar $\Omega_c^i$:
\begin{enumerate}
\item if $i\leq m$ and $j \geq m+1$, then
\begin{eqnarray*}
\widetilde{e_{i\bar j}}(z)=
\begin{cases}
\frac{1}{2}\sin^2\tau_i\bar{b_i}b_j^i & z \in \Omega_{c_1}^i\\
(\frac{1}{2}\sin^2\tau_i\bar{b_i}b_j^i)\eta(\log r_i) &
z \in \Omega_c^i \text{ and }
c_1<r_i<c\\
(\frac{1}{2}\sin^2\tau_i\bar{b_i}b_j^i)\eta(\log\rho_i-\log r_i) &
z \in \Omega_c^i \text{ and }
c^{-1}\rho_i<r_i<c_1^{-1}\rho_i\\
0 & z \in X\setminus\Omega_c^i\\
\end{cases}
\end{eqnarray*}
\item if $i,j\leq m$ and $i\ne j$, then
\begin{eqnarray*}
\widetilde{e_{i\bar j}}(z)=
\begin{cases}
\frac{1}{2}\sin^2\tau_i\bar{b_i}b_j^i & z \in \Omega_{c_1}^i\\
(\frac{1}{2}\sin^2\tau_i\bar{b_i}b_j^i)\eta(\log r_i) &
z \in \Omega_c^i \text{ and }
c_1<r_i<c\\
(\frac{1}{2}\sin^2\tau_i\bar{b_i}b_j^i)\eta(\log\rho_i-\log r_i) &
z \in \Omega_c^i \text{ and }
c^{-1}\rho_i<r_i<c_1^{-1}\rho_i\\
\frac{1}{2}\sin^2\tau_j\bar{b_i^j}b_j & z \in \Omega_{c_1}^j\\
(\frac{1}{2}\sin^2\tau_i\bar{b_i^j}b_j)\eta(\log r_j) &
z \in \Omega_c^j \text{ and }
c_1<r_j<c\\
(\frac{1}{2}\sin^2\tau_i\bar{b_i^j}b_j)\eta(\log\rho_j-\log r_j) &
z \in \Omega_c^j \text{ and }
c^{-1}\rho_j<r_j<c_1^{-1}\rho_j\\
0 & z \in X\setminus(\Omega_c^i\cup\Omega_c^j)\\
\end{cases}
\end{eqnarray*}
\item if $i\leq m$, then
\begin{eqnarray*}
\widetilde{e_{i\bar i}}(z)=
\begin{cases}
\frac{1}{2}\sin^2\tau_i |b_i|^2 & z \in \Omega_{c_1}^i\\
(\frac{1}{2}\sin^2\tau_i |b_i|^2)\eta(\log r_i) &
z \in \Omega_c^i \text{ and }
c_1<r_i<c\\
(\frac{1}{2}\sin^2\tau_i |b_i|^2)\eta(\log\rho_i-\log r_i) &
z \in \Omega_c^i \text{ and }
c^{-1}\rho_i<r_i<c_1^{-1}\rho_i\\
0 & z \in X\setminus\Omega_c^i\\
\end{cases}
\end{eqnarray*}
\end{enumerate}
Also, let $\widetilde{f_{i\bar j}}=(\Box+1)\widetilde{e_{i\bar j}}$.
It is clear that the supports of these approximation functions are
contained in the corresponding collars. We have the following estimates:
\begin{lemma}\label{etildesti}
Let $\widetilde{e_{i\bar j}}$ be the functions constructed above.
Then
\begin{enumerate}
\item $e_{i\bar i}=\widetilde{e_{i\bar i}}
+O\big (\frac{u_i^4}{|t_i|^2}\big )$ if $i \leq m$;
\item $e_{i\bar j}=\widetilde{e_{i\bar j}}
+O\big (\frac{u_i^3 u_j^3}{|t_i t_j|}\big )$ if $i,j \leq m$ and
$i\ne j$;
\item $e_{i\bar j}=\widetilde{e_{i\bar j}}
+O\big (\frac{u_i^3}{|t_i|}\big )$ if $i\leq m$ and $j\geq m+1$;
\item $\Vert e_{i\bar j}\Vert_0=O(1)$ if $i,j\geq m+1$.
\end{enumerate}
\end{lemma}

Now we use the approximation functions $\widetilde e_{i\bar j}$ in
the formulas \eqref{200}, \eqref{finalcurv} and
\eqref{finalpercurv}. The following theorems were proved in
\cite{lsy1} and \cite{lsy2}. We first have the asymptotic estimate
of the Ricci metric:
 \begin{theorem}\label{ricciest}
Let $(t,s)$ be the pinching coordinates. Then we have
\begin{enumerate}
\item $\tau_{i\bar i}=\frac{3}{4\pi^2}\frac{u_i^2}{|t_i|^2}
(1+O(u_0))$ and $\tau^{i\bar i}
=\frac{4\pi^2}{3}\frac{|t_i|^2}{u_i^2}
(1+O(u_0))$ if $i \leq m$;
\item $\tau_{i\bar j}=O\bigg (\frac{u_i^2u_j^2}{|t_it_j|}(u_i+u_j)\bigg )$ and
$\tau^{i\bar j}=O(|t_it_j|)$ if $i,j \leq m$ and $i\ne j$;
\item $\tau_{i\bar j}=O\big (\frac{u_i^2}{|t_i|}\big )$ and
$\tau^{i\bar j}=O(|t_i|)$ if $i\leq m$ and $j\geq m+1$;
\item $\tau_{i\bar j}=O(1)$ if $i,j \geq m+1$.
\end{enumerate}
\end{theorem}

By the asymptotics of the Ricci metric in the above theorem, we have
\begin{cor}\label{asypoin}
There is a constant $C>0$ such that
\[
C^{-1}\omega_{_{P}}\leq \omega_\tau\leq\omega_{_{P}}.
\]
\end{cor}

Next we estimate the holomorphic sectional curvature of the Ricci
metric:
\begin{theorem}\label{mainholo}
Let $X_0\in\bar{\mathcal{M}_g}\setminus\mathcal{M}_g$ be a
codimension $m$ point and let
$(t_1,\cdots,t_m,s_{m+1},\cdots,s_n)$ be the pinching coordinates
at $X_0$ where $t_1,\cdots,t_m$ correspond to the degeneration
directions. Then the holomorphic sectional curvature is negative
in the degeneration directions and is bounded in the
non-degeneration directions. Precisely, there is a $\delta>0$ such that if
$|(t,s)|<\delta$, then
\begin{eqnarray}\label{important100}
\widetilde R_{i\bar ii\bar i}=
\frac{3u_i^4}{8\pi^4|t_i|^4}(1+O(u_0)) >0
\end{eqnarray}
if $i\leq m$ and
\begin{eqnarray}\label{important200}
\widetilde R_{i\bar ii\bar i}=O(1)
\end{eqnarray}
if $i\geq m+1$.

Furthermore, on $\mathcal M_g$, the holomorphic sectional curvature,
the bisectional curvature and the Ricci curvature of the Ricci metric
are bounded from above and below.
\end{theorem}

This theorem was proved in \cite{lsy1} by using the formula
\eqref{finalcurv} and estimating error terms. However, the
holomorphic sectional curvature of the Ricci metric is not always
negative. We need to introduce and study the perturbed Ricci metric.
We have
\begin{theorem}\label{perholocurv}
For suitable choice of positive constant $C$, the perturbed Ricci metric
$\widetilde\tau_{i\bar j}=\tau_{i\bar j}+Ch_{i\bar j}$ is complete and
comparable with the asymptotic Poincar\'e metric.
Its bisectional curvature is bounded. Furthermore, its holomorphic
sectional curvature and Ricci curvature are bounded from
above and below by negative constants.
\end{theorem}
\begin{rem}
The perturbed Ricci metric is the first complete K\"ahler metric on
the moduli space with bounded curvature and negatively pinched
holomorphic sectional curvature and Ricci curvature.
\end{rem}

By using the minimal surface theory and Bers' embedding theorem,
we have also proved the following theorem in \cite{lsy2}:
\begin{theorem}
The moduli space equipped with either the Ricci metric or the
perturbed Ricci metric has finite volume. The Teichm\"uller space
equipped with either of these metrics has bounded geometry.
\end{theorem}

\section{The Equivalence of the Complete Metrics}

In this section we describe our arguments that all of the complete
metrics on the Teichm\"uller space and moduli space discussed
above are equivalent. With the good understanding of the Ricci and
the perturbed Ricci metrics, the results of this section are quite
easy consequences of Yau's Schwarz lemma and also the basic
definitions of these metrics. We first give the definition of
equivalence of metrics:
\begin{definition}
Two K\"ahler metrics $g_1$ and $g_2$ on a manifold $X$ are equivalent
or two norms
$\Vert \cdot\Vert_1$ and $\Vert \cdot\Vert_2$ on the tangent bundle of
$X$ are equivalent if
there is a constant $C>0$ such that
\[
C^{-1}g_1\leq g_2\leq Cg_1
\]
or
\[
C^{-1} \Vert \cdot\Vert_1\leq \Vert \cdot\Vert_2\leq C \Vert
\cdot\Vert_1.
\]
We denote this by $g_1\sim g_2$ or $\Vert \cdot\Vert_1\sim\Vert
\cdot\Vert_2$.
\end{definition}

The main result of this section we want to discuss is the following
theorem proved in \cite{lsy1} and \cite{lsy2}:
\begin{theorem}
On the moduli space $\mathcal M_g$ $(g\geq 2)$, the Teichm\"uller metric
$\Vert \cdot \Vert_T$, the Carath\'eodory metric $\Vert\cdot\Vert_C$,
the Kobayashi metric $\Vert \cdot \Vert_K$,
the K\"ahler-Einstein metric $\omega_{_{KE}}$ ,
the induced Bergman metric $\omega_{_{B}}$ , the McMullen metric
$\omega_{_{M}}$, the asymptotic Poincar\'e metric $\omega_{_{P}}$ the
Ricci metric $\omega_{\tau}$ and the perturbed Ricci metric
$\omega_{\tilde\tau}$ are equivalent. Namely
\[
\omega_{_{KE}}\sim \omega_{\tilde\tau}\sim \omega_{\tau}\sim
\omega_{_{P}}\sim \omega_{_{B}} \sim \omega_{_{M}}
\]
and
\[
\Vert \cdot \Vert_K =\Vert \cdot \Vert_T \sim  \Vert \cdot
\Vert_C\sim \Vert \cdot \Vert_{_{M}}.
\]
\end{theorem}

As corollary we proved the following conjecture of Yau made in the
early 80s \cite{yau2}, \cite{yau4}:

\begin{theorem}
The K\"ahler-Einstein metric is equivalent to the Teichm\"uller
metric on the moduli space: $\Vert \cdot \Vert_{KE}\sim\Vert \cdot
\Vert_T.$
\end{theorem}

Another corollary was also conjectured by Yau as one of his 120
famous problems \cite{yau2}, \cite{yau4}:

\begin{theorem}
The K\"ahler-Einstein metric is equivalent to the Bergman metric
on the Teichm\"uller space: $\omega_{_{KE}}\sim\omega_{_{B}}.$
\end{theorem}

Now we briefly describe the idea of proving the comparison theorem.
To compare two complete metrics on a noncompact manifold, we need to
write down their asymptotic behavior and compare near infinity.
However, if one can not find the asymptotics of these metrics, the
only tool we have is the following Yau's Schwarz lemma \cite{yau1}:
\begin{theorem}
Let $f:(M^m,g)\to (N^n,h)$ be a holomorphic map between K\"ahler
manifolds where $M$ is complete and $Ric(g)\geq -c\,g$ with $c\geq
0$.
\begin{enumerate}
\item If the holomorphic
sectional curvature of $N$ is bounded above by a negative
constant, then $f^\ast h\leq \tilde c\, g$ for some constant
$\tilde c$.
\item If $m=n$ and the Ricci curvature of $N$ is bounded above by a
negative constant, then $f^\ast\omega_h^n\leq \tilde
c\,\omega_g^n$ for some constant $\tilde c$.
\end{enumerate}
\end{theorem}

We briefly describe the proof of the comparison theorem by using
Yau's Schwarz lemma and the curvature computations and estimates.

{\bf Sketch of proof.} To use this result, we take $M=N=\M_g$ and
let $f$ be the identity map. We know the perturbed Ricci metric is
obtained by adding a positive K\"ahler metric to the Ricci metric.
Thus it is bounded from below by the Ricci metric.

Consider the identity map
\[
id:(\mathcal M_g,\omega_{\tau})\to (\mathcal M_g,\omega_{_{WP}}).
\]
Yau's Schwarz Lemma implies $\omega_{_{WP}}\leq C_0\omega_\tau$.
So
\[
\omega_\tau\leq \omega_{\tilde\tau}=\omega_\tau+C\omega_{_{WP}}\leq
(CC_0+1)\omega_\tau.
\]
Thus $\omega_\tau\sim\omega_{\widetilde\tau}$.

To control the K\"ahler-Einstein metric, we consider
\[
id:(\mathcal M_g,\omega_{_{KE}})\to (\mathcal
M_g,\omega_{\tilde\tau})
\]
and
\[
id:(\mathcal M_g,\omega_{\tilde\tau})\to (\mathcal
M_g,\omega_{_{KE}}).
\]
Yau's Schwarz Lemma implies
\[
\omega_{\tilde\tau}\leq C_0 \omega_{_{KE}}
\]
and
\[
\omega_{_{KE}}^n\leq C_0 \omega_{\tilde\tau}^n.
\]
The equivalence follows from linear algebra.

Thus by Corollary \ref{asypoin} we have
\[
\omega_{_{KE}}\sim \omega_{\tilde\tau}\sim \omega_\tau
\sim\omega_{_{P}}.
\]

By using similar method we have $\omega_\tau\leq C\omega_{_{M}}$. To
show the other side of the inequality, we have to analyze the
asymptotic behavior of the geodesic length functions. We showed in
\cite{lsy1} that
\[
\omega_\tau\sim \omega_{_{M}}.
\]
Thus by the work of McMullen \cite{mc} we have
\[
\omega_\tau\sim \omega_{_{M}}\sim \Vert\cdot\Vert_T.
\]

The work of Royden showed that the Teichm\"uller metric coincides
with the Kobayashi metric. Thus we need to show that the
Carath\'eodory metric and the Bergman metric are comparable with the
Kobayashi metric. This was done in \cite{lsy2} by using Bers'
embedding theorem. The idea is as follows:

By the Bers' embedding theorem, for each point $p\in \mathcal T_g$,
there is a map $f_p:\mathcal T_g\to \mathbb C^n$ such that $f_p(p)=0$
and
\[
B_2\subset f_p(\mathcal T_g)\subset B_6
\]
where $B_r$ is the open ball in $\mathbb C^n$ centered at $0$ with
radius $r$. Since both Carath\'eodory metric and Kobayashi metric have
restriction property and can be computed explicitly on balls, we can
use these metrics defined on $B_2$ and $B_6$ to pinch these metrics on
the Teichm\"uller space. We can also use this method to estimate peak
sections of the Teichm\"uller space at the point $p$. A careful
analysis shows
\[
\Vert\cdot\Vert_C\sim\Vert\cdot\Vert_K\sim \omega_{_{B}}.
\]
The argument is quite easy. Please see \cite{lsy2} for details.

\qed

\section{Bounded Geometry of the K\"ahler-Einstein Metric}
The comparison theorem gives us some control on the
K\"ahler-Einstein Metric. Especially we know that it has Poincar\'e
growth near the boundary of the moduli space and is equivalent to
the Ricci metric which has bounded geometry. In this section we
sketch our proof that the K\"ahler-Einstein metric also has bounded
geometry. Precisely we have
\begin{theorem}
The curvature of the K\"ahler-Einstein metric and all of its
covariant derivatives are uniformly bounded on the Teichm\"uller
spaces, and its injectivity radius has lower bound.
\end{theorem}

Now we briefly describe the proof. Please see \cite{lsy2} for
details.

{\bf Sketch of proof.} We follow Yau's argument in \cite{yau3}. The
first step is to perturb the Ricci metric using K\"ahler-Ricci flow
\[
\begin{cases}
\frac{\partial g_{i\bar j}}{\partial t}=-(R_{i\bar j}+g_{i\bar j})\\
g(0)=\tau
\end{cases}
\]
to avoid complicated computations of the covariant derivatives
of the curvature of the Ricci metric.

For $t>0$ small, let $h=g(t)$ and let $g$ be the K\"ahler-Einstein
metric. We have
\begin{enumerate}
\item $h$ is equivalent to the initial metric $\tau$ and thus is
equivalent to the K\"ahler-Einstein metric.
\item The curvature and its covariant derivatives of $h$ are bounded.
\end{enumerate}

Then we consider the Monge-Amper\'e equation
\[
\log\det(h_{i\bar j}+u_{i\bar j})-\log\det(h_{i\bar j})=u+F
\]
where $\partial\bar\partial u=\omega_{g}-\omega_h$ and
$\partial\bar\partial F=Ric(h)+\omega_h$.

The curvature of $P_{i\bar jk\bar l}$ of the K\"ahler-Einstein
metric is given by
\[
P_{i\bar jk\bar l}=R_{i\bar jk\bar l}+u_{p\bar j}h^{p\bar q}
R_{i\bar qk\bar l}+u_{;i\bar jk\bar l}-g^{p\bar q}u_{;i\bar qk}
u_{;\bar jp\bar l}.
\]

The comparison theorem implies $\partial\bar\partial u$ has
$C^0$-bound and the strong bounded geometry of $h$ implies
$\partial\bar\partial F$ has $C^k$-bound for $k\geq 0$. Also, the
equivalence of $h$ and $g$ implies $u+F$ is bounded.

So we need the $C^k$-bound of $\partial\bar\partial u$ for $k\geq
1$. Let
\[
S=g^{i\bar j}g^{k\bar l}g^{p\bar q}
u_{;i\bar qk}u_{;\bar jp\bar l}
\]
and
\begin{align*}
\begin{split}
V=& g^{i\bar j}g^{k\bar l}g^{p\bar q}g^{m\bar n}\left (
u_{;i\bar qk\bar n}u_{;\bar jp\bar lm}+
u_{;i\bar nkp}u_{;\bar jm\bar l\bar q}\right )
\end{split}
\end{align*}
where the covariant derivatives of $u$ were taken with respect to
the metric $h$.

Yau's $C^3$ estimate in \cite{yau3} implies $S$ is bounded. Let
$f=(S+\kappa)V$ where $\kappa$ is a large constant. The inequality
\[
\Delta^{'} f\geq C f^2+ (\text{ lower order terms })
\]
implies $f$ is bounded and thus $V$ is bounded. So the curvature
of the K\"ahler-Einstein metric are bounded. Same method can be used
to derive boundedness of higher derivatives of the curvature.

\qed

Actually we have also proved the all of these complete K\"ahler
metrics have bounded geometry, which should be useful in
understanding the geometry of the moduli and the Teichm\"uller
spaces.

\section{Application to Algebraic Geometry}
The existence of the K\"ahler-Einstein metric is closely related to
the stability of the tangent and cotangent bundle. In this section
 we review our results that the logarithmic extension of the cotangent bundle of the moduli
space is stable in the sense of Mumford. We first recall the
definition. Please see \cite{ko1} for details.
\begin{definition}
Let $E$ be a holomorphic vector bundle over a complex manifold $X$ and
let $\Phi$ be a K\"ahler class of $X$. The ($\Phi$-)degree of $E$ is
given by
\[
\deg(E)=\int_X c_1(E)\Phi^{n-1}
\]
where $n$ is the dimension of $X$. The slope of $E$ is given by the
quotient
\[
\mu(E)=\frac{\deg(E)}{\text{rank}(E)}.
\]
The bundle $E$ is Mumford ($\Phi$-)stable if for any proper coherent
subsheaf $\mathcal F\subset E$, we have
\[
\mu(\mathcal F)<\mu(E).
\]
\end{definition}

Now we describe the logarithmic cotangent bundle. Let $U$ be any
local chart of $\M_g$ near the boundary with pinching coordinates
$(t_1,\cdots,t_m,s_{m+1},\cdots,s_n)$ such that $(t_1,\cdots,t_m)$
represent the degeneration directions. Let
\begin{eqnarray*}
e_i=
\begin{cases}
\frac{dt_i}{t_i} & i\leq m;\\
ds_i & i\geq m+1.
\end{cases}
\end{eqnarray*}
The logarithmic cotangent bundle $E$ is the extension of
$T^\ast\M_g$ to $\bar\M_g$ such that on $U$, $e_1,\cdots,e_n$ is a
local holomorphic frame of $E$. One can write down the transition
functions and check that there is a unique bundle over $\bar\M_g$
satisfing the above condition.

To prove the stability of $E$, we need to specify a K\"ahler class.
It is natural to use the polarization of $E$. The main theorem of
this section is the following:
\begin{theorem}
The first Chern class $c_1(E)$ is
positive and $E$ is stable with respect to $c_1(E)$.
\end{theorem}

We briefly describe here the proof of this theorem. Please see
\cite{lsy2} for details.

{\bf Sketch of the proof.}
Since we only deal with the first Chern class, we can assume the
coherent subsheaf $\mathcal F$ is actually a subbundle $F$.

Since the K\"ahler-Einstein metric induces a singular metric
$g_{_{KE}}^\ast$ on the logarithmic extension bundle $E$, our main
job is to show that the degree and slope of $E$ and any proper
subbundle $F$ defined by the singular metric are finite and are
equal to the genuine ones. This depends on our estimates of the
K\"ahler-Einstein metric which are used to show the convergence of
the integrals defining the degrees.

More precisely we need to show the following:
\begin{enumerate}
\item As a current, $\omega_{_{KE}}$ is closed and represent the
first Chern class of $E$, that is
\[
[\omega_{_{KE}}]=c_1(\bar E).
\]
\item The singular metric $g_{_{KE}}^\ast$ on $E$ induced
by the K\"ahler-Einstein metric defines the degree of $E$
\[
\deg(E)=\int_{\M_g}\omega_{_{KE}}^n.
\]
\item The degree of any proper holomorphic sub-bundle $F$ of
$E$ can be defined by using $g_{_{KE}}^\ast\mid_F$,
\[
\text{deg}(F)=\int_{\M_g} -\partial\bar\partial\log\det
\left (g_{_{KE}}^\ast\mid_F\right )\wedge\omega_{_{KE}}^{n-1}.
\]
\end{enumerate}

These three steps were proved in \cite{lsy2} by using the Poincar\'e
growth property of the K\"ahler-Einstein metric together with a
special cut-off function. This shows that the bundle $E$ is
semi-stable.

To get the strict stability, we proceeded by contradiction. If $E$
is not stable, then $E$, thus $E\mid_{\M_g}$, split holomorphically.
This implies a finite smooth cover of the moduli space splits which
implies a finite index subgroup of the mapping class group splits as
a direct product of two subgroups. This is impossible by a
topological fact. Again, the detailed proof can be found in
\cite{lsy2}.

\qed
\section{Final Remarks}
Although significant progresses have been made in understanding the
geometry of the Teichm\"uller and the moduli spaces, there are still
many problems remain to be solved, such as the goodness of these
complete K\"ahler metrics, the computation of their $L^2$-cohomology
groups, the convergence of the Ricci flow starting from the Ricci
metric to the K\"ahler-Einstein metric, the representations of the
mapping class group on the middle dimensional $L^2$-cohomology of
these metrics, and the index theory associated to these complete
K\"ahler metrics. Also the perturbed Ricci metric is the first
complete K\"ahler metric on the moduli spaces with bounded negative
Ricci and holomorphic sectional curvature and bounded geometry, we
believe this metric must have more interesting applications. Another
question is which of these metrics are actually identical. We hope
to report on the progresses of the study of these problems on a
later occasion.

\end{document}